\documentclass[12pt]{article}
\usepackage{amssymb}
\usepackage{a4wide}
\usepackage{amsfonts} 
\usepackage{amsmath}
\usepackage{mathrsfs}
\newtheorem{Def}{Definition}[section]
\newtheorem{Lem}[Def]{Lemma}
\newtheorem{Thm}[Def]{Theorem}
\newtheorem{Prop}[Def]{Proposition}
\newcommand{\RR}{{\mathbb R}}
\newcommand{\CC}{{\mathbb C}}
\newcommand{\ZZ}{{\mathbb Z}}
\newcommand{\ph}{\varphi}
\newcommand{\phase}[1]{\mathrm{Phase}\,#1}
\newcommand{\id}{{\mathrm{id}}}
\newcommand{\Cb}[2]{C_{\text{\tiny bounded}}\!\left(#1,#2\right)}
\newcommand{\K}{{\mathcal K}}
\newcommand{\comp}{\!\circ\!}
\newcommand{\bez}{\setminus}
\newcommand{\spec}[1]{\mathrm{Sp}\,#1}
\newcommand{\Spec}[1]{\mathrm{Sp}\left(#1\right)}
\newcommand{\dplus}{\,\dot{+}\,}
\newcommand{\tens}{\otimes}
\newcommand{\Mor}[2]{\mathrm{Mor}\left(#1,#2\right)}
\newcommand{\aff}{\,\eta\,}
\newcommand{\M}[1]{M\!\left(#1\right)}
\newcommand{\Aut}[1]{\mathrm{Aut}\left(#1\right)}
\newcommand{\refeq}[1]{\text{(\ref{#1})}}
\newcommand{\proof}{{\sc Proof.\quad}}
\newcommand{\qed}{\hfill{\sc Q.E.D.}\medskip}
\newcommand{\atil}{\widetilde{a}}
\newcommand{\btil}{\widetilde{b}}
\newcommand{\Fq}{F_q}
\newcommand{\Y}{\Gamma}
\newcommand{\Yhat}{\widehat{\Y}}
\newcommand{\y}{\gamma}
\newcommand{\ybar}{\overline{\y}}
\newcommand{\yhat}{\widehat{\y}}
\newcommand{\Ybar}{\overline{\Y}}
\newcommand{\CinfY}{C_\infty\!\left(\Y\right)}
\newcommand{\Cinfbar}{C_\infty\!\left(\Ybar\right)}
\newcommand{\C}[2]{C\left(#1,#2\right)}
\newcommand{\HH}{\mathrm{H}}
\begin{document}
\title{Functional form of unitary representations of the quantum ``$az+b$'' 
group\thanks{Research partially supported by KBN grant No 2PO3A04022.}}
\author{W.~Pusz \& P.M.~So\l{}tan\medskip\\
{\small Department of Mathematical Methods in Physics,}\\
{\small Faculty of Physics, University of Warsaw,}\\
{\small Ho\.{z}a 74, PL 00-682 Warsaw, Poland}\\
{\small Wieslaw.Pusz@fuw.edu.pl, Piotr.Soltan@fuw.edu.pl}}
\date{}
\maketitle
\begin{abstract}
The formula for all unitary representations of the quantum ``$az+b$'' group for 
a real deformation parameter is given. The description involves the quantum 
exponential function introduced by Woronowicz.

\noindent
\textbf{Key words:} $C^*$-algebra, quantum group, unitary representation, 
exponential function.
\end{abstract}
\section{Introduction}\label{Intro}
In the study of unitary representations of a topological group it is 
interesting to look for a general structure theorem. The class of abelian 
locally compact groups provides a good example of a situation in which a general 
structure theorem is known. It is the SNAG (Stone-Naimark-Ambrose-Godement) 
theorem (\cite[Ch.~VI, Thm.~29]{gee}) saying that if $u$ is a strongly 
continuous unitary representation of a locally compact abelian group $\Y$ on a 
Hilbert space $H$ then there exists a spectral measure $E$ on the dual group 
$\Yhat$ such that
\begin{equation}\label{bichar}
u_\y=\int\limits_{\Yhat}\chi(\yhat,\y)\,dE(\yhat),
\end{equation}
where $\chi$ is a bicharacter on $\Y\times\Yhat$. The right hand side of 
\refeq{bichar} is a source of a functional expression for $u_\y$.
The best known example of such situation is described by Stone's theorem. In 
this case $\Y=\Yhat=\RR$, $\chi(\yhat,\y)=\exp{\left(i\y\yhat\right)}$ and 
\refeq{bichar} leads to a functional expression
\[
u_\y=\chi(\HH,\y)=e^{i\y\HH},
\]
where $\HH=\int\limits_{\RR}\yhat\,dE(\yhat)$ is a selfadjoint operator acting 
on $H$. No analogy of SNAG theorem is know for nonabelian case, nevertheless
in the context of locally compact quantum groups this type of structure 
theorem was proved e.g.~for the quantum $E(2)$ (cf.~\cite{e2}) and quantum 
``$ax+b$'' groups (\cite{row}). Following these ideas, in this report we shall 
describe the functional form of strongly continuous unitary representations of 
the quantum ``$az+b$'' group for real deformation parameter $q\in]0,1[$.

We shall be concerned with the theory of quantum groups on the $C^*$-algebra 
level. For the basic notions and notation used in our paper we refer the reader 
to \cite{unbo} and \cite{gen}. In particular $\M{A}$ will denote the multiplier 
algebra of a $C^*$-algebra $A$, by ``$\eta$'' we shall denote the affiliation 
relation in the sense of $C^*$-algebra theory and $\Mor{A}{B}$ will denote the 
set of morphisms from a $C^*$-algebra $A$ to a $C^*$-algebra $B$, 
i.e.~nondegenerate $*$-homomorphisms from $A$ to $M(B)$. Throughout the paper 
we shall be making extensive use of formulas (2.4)--(2.8) of \cite{gen}.
\section{Quantum ``$az+b$'' group}
The quantum ``$az+b$'' group for real deformation parameter $q\in]0,1[$ was 
described in \cite[Appendix A]{azb}. We shall denote this quantum group by the 
symbol $G$. Its construction begins with considering the set
\[
\Y=\bigl\{z\in\CC:|z|\in q^\ZZ\bigr\}
\]
and its closure in $\CC$ 
\[
\Ybar=\Y\cup\{0\}.
\]
With the topology and multiplication inherited from $\CC\bez\{0\}$, $\Y$ 
becomes an abelian topological group (isomorphic to $\ZZ\times S^1$). Clearly 
its Pontiagin dual is isomorphic to $\Y$ and the pairing describing the duality 
is given by the bicharacter
\begin{equation}\label{chi}
\chi\left(q^{i\ph+k},q^{i\psi+l}\right)=q^{i(l\ph+k\psi)}.
\end{equation}

\sloppy
Considering the natural action by multiplication of $\Y$ on $\Ybar$ we obtain a 
$C^*$-dynamical system $\left(\Cinfbar,\Y,\alpha\right)$, where the action 
$\alpha$ is
\begin{equation}\label{alact}
\bigl(\alpha_\y f\bigr)(\y')=f(\y\y')
\end{equation}
for all $f\in\Cinfbar$, $\y'\in\Ybar$ and $\y\in\Y$. The algebra $\Cinfbar$ is 
generated by an element $b\aff\Cinfbar$ given by
\begin{equation}\label{b}
b(\y)=\y
\end{equation}
for all $\y\in\Ybar$. The algebra of ``continuous functions vanishing at 
infinity on quantum ``$az+b$'' group'' is the $C^*$-algebra crossed product
\begin{equation}\label{crossed}
A=\Cinfbar\rtimes_\alpha\Y.
\end{equation}
Since the canonical injection $\Cinfbar\hookrightarrow\M{A}$ is a morphism 
from $\Cinfbar$ to $A$, we can regard $b$ as an element affiliated with $A$. By 
definition of a crossed product there is a strictly continuous family 
$(\lambda_\y)_{\y\in\Y}$ of unitary elements of $\M{A}$ implementing the action 
$\alpha$:
\begin{equation}\label{pl}
\lambda_\y f\lambda_{\y}^*=\alpha_\y(f)
\end{equation}
for any $f\in\Cinfbar$ considered as an element of $\M{A}$. Using the methods 
developed in \cite[Sect.~5]{azb} one can show that $\lambda_\y$ is of the form
\begin{equation}\label{plpl}
\lambda_\y=\chi(a,\y)
\end{equation}
for a normal element $a$ affiliated with $A$, such that $\spec{a}\subset\Ybar$.
Moreover $a$ is invertible and $a^{-1}\aff A$.

It follows from the definition of the action $\alpha$ (cf.~\refeq{pl} and 
\refeq{plpl}) that $\chi(a,\y)b\chi(a,\y)^*=\y b$ for $\y\in\Y$ and the 
$C^*$-algebra $A$ is generated by the unbounded affiliated elements 
$a,a^{-1}$ and $b$.

The basic notion used in our paper is that of a regular $q^2$-pair.
\begin{Def}\label{q2}
Let $H$ be a Hilbert space and let $(Y,X)$ be a pair of closed densely defined 
operators on $H$. We shall say that $(Y,X)$ is a regular $q^2$-pair if
\[
\left(\begin{array}{c@{\smallskip}}
X\text{ and }Y\text{ are normal},\\
\spec{X},\;\spec{Y}\subset\Ybar,\\
\ker{X}=\{0\}\text{ and}\\
\chi(X,\y)Y\chi(X,\y)^*=\y Y\\
\text{for all }\y\in\Y.
\end{array}\right)
\]
\end{Def}

Let us remark that relations of Definition \ref{q2} give the precise meaning to 
the relations of the form
\[
XY=q^2YX,\qquad XY^*=Y^*X
\]
for a pair of normal operators $(Y,X)$.

The pair $(b,a)$ of elements affiliated with the $C^*$-algebra $A$ is a regular 
$q^2$-pair in the sense that for any nondegenerate representation $\pi$ of $A$ 
on a Hilbert space $H$, the pair $(\pi(b),\pi(a))$ is a regular $q^2$-pair 
acting on $H$. We shall also use the notion of a regular $q^2$-pair in the 
context of any $C^*$-algebra $B$.

Let us note that the $C^*$-algebra $A$ defined by \refeq{crossed} is a 
universal $C^*$-algebra generated by a regular $q^2$-pair in the following 
sense (comp.~\cite[Prop.~4.2]{azb}):

\begin{Prop}\label{uni}
Let $B$ be a $C^*$-algebra and $(b_0,a_0)$ a regular $q^2$-pair of elements 
affiliated with $B$. Then there exists a unique morphism 
$\ph\in\Mor{A}{B}$ such that
\[
\begin{array}{r@{\;=\;}l@{\smallskip}}
\ph(a)&a_0,\\
\ph(b)&b_0.
\end{array}
\]
\end{Prop}

It turns out that the operator $a\tens b+b\tens I$ is closeable and (denoting 
its closure by $a\tens b\dplus b\tens I$ the pair 
$(a\tens b\dplus b\tens I,a\tens a)$ is a regular $q^2$-pair 
(cf.~Proposition \ref{Qexp}) of elements affiliated with $A\tens A$. 
Therefore there exists a unique morphism $\Delta\in\Mor{A}{A\tens A}$ such that
\begin{equation}\label{DelG}
\begin{array}{r@{\;=\;}l@{\smallskip}}
\Delta(a)&a\tens a,\\
\Delta(b)&a\tens b\dplus b\tens I.
\end{array}
\end{equation}
Moreover $\Delta$ is coassociative and encodes the group structure of the 
quantum ``$az+b$'' group, briefly $G=(A,\Delta)$.
\section{Quantum exponential function}
In \cite{opeq} S.L.~Woronowicz introduced the quantum exponential function 
$\Fq$ defined on $\Ybar$. It is given by the formula
\[
\Fq(\y)=\prod_{k=0}^\infty\frac{1+q^{2k}\ybar}{1+q^{2k}\y}
\]
for $\y\in\Ybar\bez\{-1,-q^{-2},-q^{-4},\ldots\}$ and $\Fq(\y)=-1$ for
$\y\in\{-1,-q^{-2},-q^{-4},\ldots\}$. Thus defined, $\Fq$ is a 
continuous function $\Ybar\to S^1$. Moreover $\Fq(0)=1$.

The most important property of the quantum exponential function is the 
one contained in the following proposition:

\begin{Prop}[{\cite[Thm.~3.1]{opeq}}]\label{Qexp}
Let $H$ be a Hilbert space and let $(Y,X)$ be a regular $q^2$-pair acting on 
$H$. Then the sum $X+Y$ is a densely defined closeable operator and its closure 
$X\dplus Y$ is a normal operator with $\Spec{X\dplus Y}\subset\Ybar$. Moreover
\begin{equation}\label{EXP}
\Fq(X\dplus Y)=\Fq(Y)\Fq(X).
\end{equation}
\end{Prop}

The last statement in Proposition \ref{Qexp} justifies the name ``quantum 
exponential function''. 

Let us remark that formula \refeq{EXP} holds for more general $q^2$-pairs, 
without the assumption that $\ker{X}=\{0\}$ (cf.~\cite[Section 2]{gl}). 
Moreover $\Fq$ is the only solution of this type of functional equation in 
a more general sense. To formulate the corresponding result let $H$ be a Hilbert 
space and $f\colon\Ybar\ni\y\mapsto f(\y)\in B(H)$ be a bounded measureable 
mapping. For a normal operator $Y$ acting on a Hilbert space $K$ such that 
$\spec{Y}\subset\Ybar$ we set
\[
f(Y)=\int\limits_{\Ybar}f(\y)\tens dE_Y(\y),
\]
where $dE_Y(\y)$ is the spectral measure of $Y$. Clearly $f(Y)\in B(H\tens K)$.
\begin{Thm}[{\cite[Thm.~4.2]{opeq}}]\label{ExpT}
Let $H$ and $K$ be Hilbert spaces and let $(Y,X)$ be a regular $q^2$-pair 
acting on $K$. Let 
\[
f\colon\Ybar\ni\y\longmapsto f(\y)\in B(H)
\]
be a borel mapping such that $f(\y)$ is unitary for almost all $\y\in\Ybar$. 
Assume that
\begin{equation}\label{dd}
f(X\dplus Y)=f(Y)f(X).
\end{equation}
Then there exists a normal operator $Z$ on $H$ such that $\spec{Z}\subset\Ybar$ 
and 
\[
f(\y)=\Fq(Z\y)
\]
for almost all $\y\in\Ybar$.
\end{Thm}

Note that any borel solution of the functional equation \refeq{dd} is in fact a 
continuous one. This result is crucial for functional description of unitary 
representations of $G$ given in the next section.
\section{Structure of unitary representations of the quantum ``$az+b$'' group}
In this section we investigate strongly continuous unitary representations of 
$G$.
\begin{Def}\label{URep}
A strongly continuous unitary representation of $G$ on a Hilbert space $H$ is a 
unitary element 
\[
U\in\M{\K(H)\tens A}
\]
such that
\begin{equation}\label{Urep}
(\id\tens\Delta)U=U_{12}U_{13}
\end{equation}
(where we used the leg numbering notation).
\end{Def}

In what follows we shall abbreviate ``strongly continuous unitary 
representations'' to ``unitary representations''.
The main result of the paper is contained in the following theorem:
\begin{Thm}\label{glowne}
Let $U$ be a unitary representation of the quantum ``$az+b$'' group on a Hilbert 
space $H$. Then there exists a unique regular $q^2$-pair 
$(\btil,\atil)$ acting on $H$ such that
\begin{equation}\label{U}
U=\Fq(\btil\tens b)\chi(\atil\tens I,I\tens a).
\end{equation}
Conversely, for any regular $q^2$-pair $(\atil,\btil)$ acting on a Hilbert 
space $H$, the operator $U$ defined by formula \refeq{U} is a unitary 
representation of $G$.
\end{Thm}
\proof
Let
\[
\begin{array}{r@{\;=\;}l@{\smallskip}}
a_0(\y)&\y,\\
b_0(\y)&0
\end{array}
\]
for all $\y\in\Y$. Then $a_0,b_0\aff\CinfY$ and if we represent $\CinfY$
on $L^2(\Y)$ by multiplication operators, $(b_0,a_0)$ becomes a regular 
$q^2$-pair on $L^2(\Y)$. By the universal property of $A$ (Proposition 
\ref{uni}) there exists a 
$\ph\in\Mor{A}{\CinfY}=C\left(\Y,\Mor{A}{\CC}\right)$ such that
\[
\begin{array}{r@{\;=\;}l@{\smallskip}}
\ph(a)&a_0,\\
\ph(b)&b_0.
\end{array}
\]
Thus $\ph$ is a continuous family $\ph=(\ph_\y)_{\y\in\Y}$ where the $\ph_\y$ 
are multiplicative functionals on $A$:
\[
\ph_\y(x)=\bigl(\ph(x)\bigr)(\y).
\]
Moreover the map $\Y\ni\y\mapsto\ph_\y\in A^*$ is a homomorphism in the sense 
that
\begin{equation}\label{splotph}
\ph_{\y_1}*\ph_{\y_2}=\ph_{\y_1\y_2},
\end{equation}
where the convolution
\[
\ph_{\y_1}*\ph_{\y_2}=(\ph_{\y_1}\tens\ph_{\y_2})\comp\Delta.
\]
Let 
\[
\phi(x)=(\id\tens\ph)\Delta(x)
\]
for $x\in A$. Then $\phi\in\Mor{A}{A\tens\CinfY}=\C{\Y}{\Mor{A}{A}}$. Again we 
identify $\phi$ with a continuous family $\phi=(\phi_\y)_{\y\in\Y}$ and it is 
easy to see that for all $\y\in\Y$ the maps $\phi_\y$ are automorphisms of $A$:
\[
\phi_\y\in\Aut{A}.
\]
In other words
\[
\phi\in\C{\Y}{\Aut{A}}.
\]
It follows from \refeq{splotph} that $(\phi_\y)_{\y\in\Y}$ is a continuous 
group of automorphisms of $A$. It is also easy to check that 
\[
\begin{array}{r@{\;=\;}l@{\smallskip}}
\phi_\y(a)&\y a,\\
\phi_\y(b)&b
\end{array}
\]
for all $\y\in\Y$. In other words the action 
$\Yhat=\Y\ni\y\mapsto\phi_\y\in\Aut{A}$ is the dual action to the action 
$\alpha$ of $\Y$ on $\Cinfbar$ (cf.~\refeq{alact} and \refeq{crossed}).

Define 
\[
u=(\id\tens\ph)U\in\M{\K(H)\tens\CinfY}=\Cb{\Y}{B(H)}.
\]
We can thus view $u$ as a continuous family $(u_\y)_{\y\in\Y}$ of unitary 
elements in $B(H)$.

Since by \refeq{Urep} and \refeq{splotph}
\[
\begin{array}{r@{\;=\;}l@{\smallskip}}
u_{\y_1}u_{\y_2}&\bigl((\id\tens\ph_{\y_1})U\bigr)
\bigl((\id\tens\ph_{\y_2})U\bigr)\\
&(\id\tens\ph_{\y_1}\tens\ph_{\y_2})U_{12}U_{13}\\
&(\id\tens\ph_{\y_1}*\ph_{\y_2})U=u_{\y_1\y_2},
\end{array}
\]
$u$ is a strongly continuous representation of $\Y$ in the Hilbert space $H$. 
By SNAG theorem (cf.~\cite[Ch.~VI, Thm.~29]{gee} and Section \ref{Intro}) 
\[
u_\y=\chi(\atil,\y)
\]
where $\atil$ is a normal operator on $H$ such that $\spec{\atil}\subset\Ybar$ 
and $\ker{\atil}=\{0\}$.

Let
\[
V=\chi(\atil\tens I,I\tens a)\in\M{\K(H)\tens A}.
\]
We have
\[
(\id\tens\ph)V=u.
\]

Using \refeq{Urep} we obtain
\begin{equation}\label{phU}
\begin{array}{r@{\;=\;}l@{\smallskip}}
(\id\tens\phi_\y)U&(\id\tens\id\tens\ph_{\y})\Delta(U)\\
&U\bigl((\id\tens\ph_{\y})U\bigr)\\
&U(u_\y\tens I)=U\bigl(\chi(\atil,\y)\tens I\bigr)
\end{array}
\end{equation}
and by \refeq{DelG} and the definition of $\ph_\y$ we have
\[
\begin{array}{r@{\;=\;}l@{\smallskip}}
(\id\tens\phi_\y)V&
(\id\tens\id\tens\ph_{\y})(\id\tens\Delta)\chi(\atil\tens I,I\tens a)\\
&\chi(\atil\tens I,I\tens\y a)=V\bigl(\chi(\atil,\y)\tens I\bigr).
\end{array}
\]
Define 
\[
W=UV^*.
\]
It follows that for all $\y\in\Y$
\begin{equation}\label{inv}
(\id\tens\phi_\y)W=W.
\end{equation}
At this point one expects that 
\begin{equation}\label{nal}
W\in\M{\K(H)\tens\Cinfbar}.
\end{equation}
It is known that $\M{\K(H)\tens\Cinfbar}=\Cb{\Ybar}{B(H)}$, therefore $W=f(b)$ 
where $f\in\Cb{\Ybar}{B(H)}$ and $f(z)$ is unitary for any $z\in\Ybar$. 
Unfortunately in the context of $C^*$-algebra crossed products the invariance 
condition \refeq{inv} is not sufficient to support \refeq{nal}. In addition to 
\refeq{inv} one needs to know that 
$\Y\ni\y\mapsto(I\tens\lambda_\y)W(I\tens\lambda_\y)^*X$ is norm continuous for 
any $X\in\K(H)\tens\Cinfbar\subset\M{\K(H)\tens A}$ 
(cf.~\cite[Proposition 7.8.9]{ped}). However we have no argument to justify 
this. On the other hand for $W^*$-dynamical systems the additional condition is 
not relevant. Then we only have borel measureability of the corresponding 
function $f$. Nevertheless this fact combined with a functional equation for 
$f$ (cf.~\refeq{dd}) will imply continuity of $f$. Therefore we are able to 
justify \refeq{nal} only {\em a posteriori}.

Let us now extend the $C^*$-dynamical system 
$\left(\K(H)\tens\Cinfbar,\Y,\id\tens\alpha\right)$ to a $W^*$-dynamical system
$\left(B(H)\tens L^\infty\left(\Ybar\right),\Y,\id\tens\alpha\right)$. We may
assume that $\Cinfbar$ is faithfully represented in a Hilbert space.
Then 
\[
W\in\M{\K(H)\tens A}\subset B(H)\tens A''=
\left(B(H)\tens L^\infty\left(\Ybar\right)\right)\rtimes_{\id\tens\alpha}\Y
\]
and $W$ is fixed under the action $\y\mapsto(\id\tens\phi_\y)$ which clearly 
is dual to the action $\y\mapsto(\id\tens\alpha_\y)$. By 
\cite[Theorem 7.10.4]{ped} the element $W$ belongs to the von Neuman algebra 
\[
B(H)\tens L^\infty\left(\Ybar\right)
\]
which means that $W=f(b)$ where $f\colon\Ybar\to B(H)$ is a unitary 
operator-valued borel function. Now we shall show that $f$ satisfies the 
functional equation \refeq{dd}.

Since $U=WV=f(b)\chi(\atil\tens I,I\tens a)$ we have by \refeq{DelG}
\begin{equation}\label{DelU}
(\id\tens\Delta)U=f(a\tens b\dplus b\tens I)
\chi(\atil\tens I\tens I,I\tens a\tens a).
\end{equation}
Therefore
\[
\begin{array}{r@{\;=\;}l@{\smallskip}}
(\id\tens\phi_\y\tens\id)U_{12}U_{13}
&(\id\tens\ph_\y\tens\id)(\id\tens\Delta)U\\
&f(\y b)\chi(\atil\tens I,I\tens\y a)\\
&f(\y b)\bigl(\chi(\atil,\y)\tens I\bigr)\chi(\atil\tens I,I\tens a).
\end{array}
\]
On the other hand (cf.~\refeq{phU})
\[
(\id\tens\ph_\y\tens\id)U_{12}U_{13}=\bigl(\chi(\atil,\y)\tens I\bigr)U
\]
and thus
\begin{equation}\label{tozs}
\begin{array}{r@{\;=\;}l@{\smallskip}}
U&\bigl(\chi(\atil,\y)\tens I\bigr)^*
f(\y b)\bigl(\chi(\atil,\y)\tens I\bigr)\chi(\atil\tens I,I\tens a)\\
&\chi(\atil\tens I,\y I\tens I)^*
f(\y b)\chi(\atil\tens I,\y I\tens I)\chi(\atil\tens I,I\tens a)
\end{array}
\end{equation}
for all $\y\in\Y$. But $U$ does not depend on $\y$, so integrating both sides 
of \refeq{tozs} over $\y$ with respect to the spectral measure of $a$ we obtain
\[
\begin{array}{r@{\;=\;}l@{\smallskip}}
U_{13}&\chi(\atil\tens I\tens I,I\tens a\tens I)^*f(a\tens b)
\chi(\atil\tens I\tens I,I\tens a\tens I)
\chi(\atil\tens I\tens I,I\tens I\tens a)\\
&\chi(\atil\tens I\tens I,I\tens a\tens I)^*f(a\tens b)
\chi(\atil\tens I\tens I,I\tens a\tens a).
\end{array}
\]
Now since 
\[
U_{12}=f(b\tens I)\chi(\atil\tens I\tens I,I\tens a\tens I)
\]
we obtain
\begin{equation}\label{UU}
U_{12}U_{13}=f(b\tens I)f(a\tens b)\chi(\atil\tens I\tens I,I\tens a\tens a).
\end{equation}
Comparing \refeq{UU} with \refeq{DelU} we obtain a functional expression
\begin{equation}\label{exponf}
f(a\tens b\dplus b\tens I)=f(b\tens I)f(a\tens b).
\end{equation}
Denote $X=b\tens I$, $Y=a\tens b$. It is easy to check that $(Y,X)$ is a 
regular $q^2$-pair. Therefore by \refeq{exponf} and Theorem \ref{ExpT} there 
exists a normal operator $\btil$ acting on $H$ with $\spec{\btil}\subset\Ybar$ 
such that
\[
f(\y)=\Fq(\btil\y).
\]
Consequently
\[
f(b)=\Fq(\btil\tens b)
\]
and
\begin{equation}\label{postac}
U=\Fq(\btil\tens b)\chi(\atil\tens I,I\tens a).
\end{equation}

So far we know that $\atil$ and $\btil$ are normal operators on $H$ such that 
$\spec{\atil},\;\spec{\btil}\subset\Ybar$ and $\ker{\atil}=\{0\}$. To end the 
proof of the existence part of Theorem \ref{glowne} cf.~Definition \ref{q2}) 
we need to show that 
\begin{equation}\label{weyl}
\chi(\atil,\y)\btil\chi(\atil,\y)^*=\y\btil
\end{equation}
for all $\y\in\Y$.

To that end observe that inserting the information about $U$ given by 
\refeq{postac} into the identity \refeq{tozs} we obtain
\[
\bigl(\chi(\atil,\y)\tens I\bigr)\Fq(\btil\tens b)
\bigl(\chi(\atil,\y)\tens I\bigr)^*=\Fq(\y\btil\tens b)
\]
which by unitarity of $\bigl(\chi(\atil,\y)\tens I\bigr)$ means that
\begin{equation}\label{prawie}
\Fq\left(\bigl(\chi(\atil,\y)\btil\chi(\atil,\y)^*\bigr)\tens b\right)
=\Fq(\y\btil\tens b).
\end{equation}
Setting $T_1=\chi(\atil,\y)\btil\chi(\atil,\y)^*$ and $T_2=\y\btil$ we can 
rewrite \refeq{prawie} as
\begin{equation}\label{lepiej}
\Fq(T_1\tens b)=\Fq(T_2\tens b).
\end{equation}
For $z\in\Ybar$ let $\omega_z\in\Mor{\Cinfbar}{\CC}$ be given by
\[
\omega_z(b)=z.
\]
Applying $(\id\tens\omega_z)$ to both sides of \refeq{lepiej} we get
\[
\Fq(zT_1)=\Fq(zT_2)
\]
and this equality holds for any $z\in\Ybar$. Now the equality of $T_1$ and 
$T_2$ follows from the following result:
\begin{Lem}\label{iff}
Let $T_1$ and $T_2$ be normal operators acting on a Hilbert space $K$ such that 
$\spec{T_1},\;\spec{T_2}\subset\Ybar$. Then
\[
\left(\begin{array}{c}\Fq(zT_1)=\Fq(zT_2)\\\text{for all }z\in\Ybar
\end{array}\right)\Longleftrightarrow\Bigl(T_1=T_2\Bigr).
\]
\end{Lem}
(We omit the proof since it is analogous to the proof of \cite[Lemma 3.5]{gl}.)

Therefore $(\btil,\atil)$ is a regular $q^2$-pair. 

To prove uniqueness of $(\btil,\atil)$ let us apply $(\id\tens\ph_\y)$ to $U$:
\[
(\id\tens\ph_\y)U=\chi(\atil,\y).
\]
It shows that the operator $\atil$ is determined uniquely. Indeed: we have 
(cf.~\refeq{chi})
\[
\begin{array}{r@{\;=\;}l@{\smallskip}}
\chi(\y,q)&\phase{\y},\\
\chi(\y,q^{it})&|\y|^{it}
\end{array}
\]
for all $\y\in\Y$. Therefore, by functional calculus for normal operators,
\[
\begin{array}{r@{\;=\;}l@{\smallskip}}
\phase{\atil}&(\id\tens\ph_q)U,\\
|\atil|^{it}&(\id\tens\ph_{q^{it}})U
\end{array}
\]
which determines $\atil$ completely. Now the operator $\btil$ is also 
determined uniquely. In fact if
\[
U=\Fq(\btil'\tens b)\chi(\atil\tens I,I\tens a)
\]
then
\[
\Fq(\btil'\tens b)=\Fq(\btil\tens b)
\]
and the reasoning presented after \refeq{lepiej} shows that $\btil'=\btil$. 
This ends the proof of the first part of our theorem.

For the proof of the second part let $(\btil,\atil)$ be a regular $q^2$-pair 
acting on a Hilbert space $H$ and
\[
U=\Fq(\btil\tens b)\chi(\atil\tens I,I\tens a).
\]
Elements $\btil\tens b$, $\atil\tens I$ and $I\tens a$ are affiliated 
with $\K(H)\tens A$. Therefore $U\in\M{\K(H)\tens A}$ and $U$ is unitary. Now 
by \refeq{DelG}
\begin{equation}\label{this}
(\id\tens\Delta)U=\Fq(\btil\tens a\tens b\dplus\btil\tens b\tens I)
\chi(\atil\tens I\tens I,I\tens a\tens a).
\end{equation}
Since $(\btil\tens b\tens I,\btil\tens a\tens I)$ is a $q^2$-pair, we 
have (cf.~remark after Proposition \ref{Qexp})
\begin{equation}\label{this2}
\Fq(\btil\tens a\tens b\dplus\btil\tens b\tens I)
=\Fq(\btil\tens b\tens I)\Fq(\btil\tens a\tens b).
\end{equation}
Moreover by the character property of $\chi$
\[
\chi(\atil\tens I\tens I,I\tens a\tens a)=
\chi(\atil\tens I\tens I,I\tens a\tens I)
\chi(\atil\tens I\tens I,I\tens I\tens a)
\]
and $\y\btil\,\chi(\atil,\y)=\chi(\atil,\y)\,\btil$ by \refeq{weyl}. Therefore
\[
(\btil\tens a)\chi(\atil\tens I,I\tens a)=
(\chi(\atil\tens I,I\tens a)(\btil\tens I)
\]
and
\begin{equation}\label{klucz}
(\btil\tens a\tens b)\chi(\atil\tens I\tens I,I\tens a\tens I)
=\chi(\atil\tens I\tens I,I\tens a\tens I)(\btil\tens I\tens b).
\end{equation}
Now
\[
\begin{array}{l@{\smallskip}}
\Fq(\btil\tens b\tens I)\Fq(\btil\tens a\tens b)
\chi(\atil\tens I\tens I,I\tens a\tens a)\\
\quad=\Fq(\btil\tens b\tens I)\Fq(\btil\tens a\tens b)
\chi(\atil\tens I\tens I,I\tens a\tens I)
\chi(\atil\tens I\tens I,I\tens I\tens a)\\
\quad=\Fq(\btil\tens b\tens I)\chi(\atil\tens I\tens I,I\tens a\tens I)
\Fq(\btil\tens I\tens b)\chi(\atil\tens I\tens I,I\tens I\tens a)\\
\quad=U_{12}U_{13}.
\end{array}
\]
This combined with \refeq{this} and \refeq{this2} shows that $U$ is a unitary 
representation of $G$ (cf.~Definition \ref{URep}).
%
\qed
\section{Acknowledgements}
The authors wish to thank S.L.~Woronowicz for many stimulating discussions 
on the subject of quantum groups and their representations.

\end{document}